\documentclass[preprint,12pt]{elsarticle}

\usepackage[T1]{fontenc}
\usepackage[utf8]{inputenc}
\usepackage{lmodern}
\usepackage{amsmath,amssymb,amsthm}
\usepackage{booktabs}
\usepackage{microtype}
\usepackage{enumitem}
\usepackage[hidelinks]{hyperref}

\journal{J. Math. Anal. Appl.}

\newtheorem{theorem}{Theorem}[section]
\newtheorem{proposition}[theorem]{Proposition}
\newtheorem{corollary}[theorem]{Corollary}
\newtheorem*{theoremA}{Theorem A}
\newtheorem*{theoremB}{Theorem B}
\theoremstyle{remark}
\newtheorem{remark}[theorem]{Remark}

\newcommand{\R}{\mathbb{R}}
\newcommand{\Crit}{\operatorname{Crit}}

\begin{document}

\begin{frontmatter}

\title{The local-to-global minimum problem of Calvert and Vamanamurthy: a quartic counterexample in three variables and the cubic case in all dimensions}

\author[unicordoba]{Hugo Adu\'en\corref{cor1}}
\ead{haduen@correo.unicordoba.edu.co}
\cortext[cor1]{Corresponding author.}
\address[unicordoba]{Departamento de Matem\'aticas y Estad\'istica, Universidad de C\'ordoba, Carrera 6 No.\ 77-305, Monter\'ia, C\'ordoba 230002, Colombia}

\begin{abstract}
In 1980 Calvert and Vamanamurthy asked whether a polynomial whose only critical point is a local minimum must attain its global minimum there. They gave counterexamples of every degree at least five in two or more variables and settled cubics in two and three variables and quartics in two variables, leaving open quartics in three or more variables and cubics in four or more variables. We resolve both remaining cases. First, we give an explicit quartic in three variables whose unique critical point is a nondegenerate strict local minimum, although the polynomial is unbounded below; the example extends to every dimension at least three, and its critical locus is a single point even over the complex numbers. Second, we prove that in every dimension a polynomial of exact degree three with a local minimum has another critical point, so that the cubic case holds vacuously; the proof is elementary and avoids the topological degree arguments used previously in two and three variables. Together with the results of Calvert and Vamanamurthy, this gives the complete classification of the local-to-global property by degree and dimension.
\end{abstract}


\end{frontmatter}

\section{Introduction}\label{sec:intro}
Let $p:\R^n\to\R$ be a real polynomial. We consider the following local-to-global question:
\begin{quote}
If $p$ has a local minimum at its only critical point, must that point be a global minimum?
\end{quote}
Write $\Crit(p):=\{x\in\R^n:\nabla p(x)=0\}$ for the set of critical points of $p$. By translating the variables and subtracting a constant, the question takes the equivalent form
\begin{equation}\label{eq:question}
\begin{gathered}
\Crit(p)=\{0\}\quad\text{and}\quad 0\ \text{is a local minimum of}\ p\\[2pt]
\Longrightarrow\quad p(x)\ge p(0)\quad\text{for all } x\in\R^n.
\end{gathered}
\end{equation}

For one variable, the implication follows directly from elementary calculus. In several variables it can fail because a polynomial may decrease along a path leading to infinity without acquiring another finite critical point. This phenomenon is a classical source of examples in the teaching of calculus \cite{RosenholtzSmylie1985,AshSexton1985}, it is closely related to the behaviour of a polynomial at infinity \cite{Durfee1998}, and it reappears in the variational literature \cite{Bisgard2015} and in recent work on local extrema and local minima of polynomials \cite{Fernando2020,AhmadiZhang2022}.

Calvert and Vamanamurthy \cite{CalvertVamanamurthy1980} investigated this question systematically. Their Theorem~1 gives a counterexample of each degree $d\ge5$ in every dimension $n\ge2$. Their Theorem~5 proves \eqref{eq:question} for quartics in two variables, and their Theorems~4 and~6 show that a cubic in two or three variables with a local minimum has another critical point, so that \eqref{eq:question} holds for such cubics. They stated explicitly that the remaining cases, namely \emph{quartics in three or more variables} and \emph{cubics in four or more variables}, were open.

The purpose of this paper is to settle both remaining cases.

\begin{theoremA}
For every $n\ge3$ there is a polynomial of degree four on $\R^n$ whose only critical point is a nondegenerate strict local minimum, although the polynomial is unbounded below.
\end{theoremA}

\begin{theoremB}
In every dimension $n\ge1$, a polynomial of exact degree three with a local minimum has a critical point different from that minimum. Consequently, \eqref{eq:question} holds, vacuously, for cubic polynomials in every dimension.
\end{theoremB}

Theorem~A is proved in Section~\ref{sec:example} (Theorem~\ref{thm:quartic} and Corollary~\ref{cor:high-dim}) by means of the explicit quartic
\[
p(x,y,z)=y^2+(1+y)z^2+\bigl(z-x(1+y)\bigr)^2.
\]
Its quadratic part at the origin is positive definite, so the counterexample persists under the stronger hypothesis of a nondegenerate strict local minimum. Moreover, the origin is its only critical point even in $\mathbb{C}^3$ (Remark~\ref{rem:complex}): the failure of \eqref{eq:question} is not caused by complex critical points but by the behaviour of $p$ at infinity. Theorem~B is proved in Section~\ref{sec:cubic} (Theorem~\ref{thm:cubic}) by a Lagrange multiplier argument combined with Euler's identity; it contains Theorems~4 and~6 of \cite{CalvertVamanamurthy1980} as special cases and does not use topological degree.

Combining Theorems~A and~B with the results of \cite{CalvertVamanamurthy1980} gives the sharp degree--dimension classification in Section~\ref{sec:classification}. Section~\ref{sec:elementary} records elementary proofs for degrees one and two, and Section~\ref{sec:construction} explains how the quartic arises from the degree-five example of \cite{CalvertVamanamurthy1980}.

\section{A quartic counterexample}\label{sec:example}
\begin{theorem}\label{thm:quartic}
For any $a,b>0$, define $p_{a,b}:\R^3\to\R$ by
\begin{equation}\label{eq:family}
p_{a,b}(x,y,z)=a y^2+(1+y)z^2+b\bigl(z-x(1+y)\bigr)^2.
\end{equation}
Then $p_{a,b}$ has degree exactly four and has the following properties:
\begin{enumerate}[label=\upshape(\roman*),leftmargin=2.4em,itemsep=1pt,topsep=3pt]
\item the origin is a nondegenerate strict local minimum;
\item $\Crit(p_{a,b})=\{(0,0,0)\}$;
\item $p_{a,b}$ is unbounded below.
\end{enumerate}
\end{theorem}
\begin{proof}
The coefficient of $x^2y^2$ in \eqref{eq:family} is $b>0$, so the polynomial has degree four. Its quadratic part at the origin is
\[
Q(x,y,z)=a y^2+z^2+b(z-x)^2,
\]
which is positive definite. Hence $D^2p_{a,b}(0)$ is positive definite and the origin is a nondegenerate strict local minimum.

Put $h=z-x(1+y)$. Direct differentiation gives
\begin{equation}\label{eq:partials}
\begin{aligned}
\partial_xp_{a,b}&=-2b(1+y)h,\\
\partial_yp_{a,b}&=2ay+z^2-2bxh,\\
\partial_zp_{a,b}&=2(1+y)z+2bh.
\end{aligned}
\end{equation}
Suppose $\nabla p_{a,b}(x,y,z)=0$. If $y=-1$, the last equation in \eqref{eq:partials} gives $h=z=0$, whereas the middle equation gives $\partial_yp_{a,b}=-2a\ne0$. This is impossible. Thus $1+y\ne0$ and the first equation gives $h=0$, that is, $z=x(1+y)$. The last equation then becomes $2(1+y)z=0$, whence $z=0$ and $x=0$. Finally, the middle equation yields $y=0$. Conversely, the origin is critical. This proves the uniqueness assertion.

For every $t\in\R$, on the line $(x,y,z)=(t,-2,-t)$ we have
\begin{equation}\label{eq:negative}
p_{a,b}(t,-2,-t)=4a-t^2\to-\infty
\qquad\text{as }|t|\to\infty.
\end{equation}
Thus the local minimum is not global.
\end{proof}

The particularly simple choice $a=b=1$ gives
\begin{equation}\label{eq:main-example}
\boxed{\quad p(x,y,z)=y^2+(1+y)z^2+\bigl(z-x(1+y)\bigr)^2.\quad}
\end{equation}
In this case $p(3,-2,-3)=-5<p(0,0,0)=0$.

\begin{remark}\label{rem:complex}
The proof of the uniqueness assertion in Theorem~\ref{thm:quartic} uses only the field operations and the facts $a\ne0$, $b\ne0$; it never uses that the variables are real. Hence the same argument shows that the origin is the only critical point of $p_{a,b}$ in $\mathbb{C}^3$. In particular, the failure of \eqref{eq:question} for $p_{a,b}$ cannot be attributed to complex critical points; it is entirely a phenomenon at infinity, in the spirit of \cite{Durfee1998}.
\end{remark}

\begin{corollary}\label{cor:high-dim}
For every integer $n\ge3$, there exists a polynomial of degree four on $\R^n$ with exactly one critical point, which is a nondegenerate strict local minimum but not a global minimum.
\end{corollary}
\begin{proof}
For $n>3$, set
\[
P(x,y,z,w_4,\dots,w_n)=p(x,y,z)+\sum_{j=4}^n w_j^2,
\]
where $p$ is given by \eqref{eq:main-example}. Its critical-point equations force $w_4=\cdots=w_n=0$ and $\nabla p(x,y,z)=0$. Its quadratic part is positive definite, and \eqref{eq:negative} with $a=b=1$ gives a sequence along which $P\to-\infty$.
\end{proof}

\section{The cubic case in arbitrary dimension}\label{sec:cubic}
The following theorem settles the cubic case left open in \cite{CalvertVamanamurthy1980} for four or more variables. It contains Theorems~4 and~6 of \cite{CalvertVamanamurthy1980}, which treat two and three variables by means of topological degree, and its proof is elementary.

\begin{theorem}\label{thm:cubic}
Let $p:\R^n\to\R$ be a polynomial of exact degree three. If $0$ is a local minimum of $p$, then $p$ has a critical point different from $0$.
\end{theorem}
\begin{proof}
Subtracting $p(0)$ and using $\nabla p(0)=0$, write
\[
p(x)=Q(x)+C(x),
\]
where $Q$ and $C$ are homogeneous of degrees two and three and $C\not\equiv0$. Local minimality implies $Q\ge0$.

Suppose first that $Q$ is not positive definite. Choose $v\ne0$ with $Q(v)=0$. Since $p(tv)=t^3C(v)$, local minimality for both signs of $t$ gives $C(v)=0$. Thus $p(tv)=p(0)=0$ for all $t$. Every sufficiently small nonzero point $tv$ is also a local minimum (because $p\ge p(0)$ in a fixed neighborhood of $0$) and hence a nonzero critical point.

Suppose next that $Q$ is positive definite. Then the ellipsoid
\[
E=\{v\in\R^n:Q(v)=1\}
\]
is compact. Because $C$ is nonzero and odd, it attains a negative minimum at some $v\in E$, and at this point the Lagrange multiplier rule gives $\nabla C(v)=\lambda\nabla Q(v)$ for some $\lambda\in\R$. Euler's identities for homogeneous polynomials imply
\[
3C(v)=\langle\nabla C(v),v\rangle
=\lambda\langle\nabla Q(v),v\rangle=2\lambda,
\]
so $\lambda=\frac32C(v)<0$. With $r=-1/\lambda>0$, we obtain
\[
\nabla p(rv)=r\nabla Q(v)+r^2\nabla C(v)
=r(1+r\lambda)\nabla Q(v)=0.
\]
Since $rv\ne0$, this is the required second critical point.
\end{proof}

\section{Sharp classification by degree and dimension}\label{sec:classification}
For integers $n\ge1$ and $d\ge2$, let $\mathcal P(n,d)$ denote the assertion that every polynomial $p:\R^n\to\R$ of degree at most $d$ satisfying the two hypotheses in \eqref{eq:question} also satisfies its conclusion.

\begin{theorem}\label{thm:classification}
For $n\ge1$ and $d\ge2$,
\begin{equation}\label{eq:threshold}
\boxed{\;\mathcal P(n,d)\ \textup{holds}
\iff
\begin{cases}
n=1, & \textup{or}\\
n=2\ \textup{and}\ d\le4, & \textup{or}\\
n\ge3\ \textup{and}\ d\le3.
\end{cases}\;}
\end{equation}
\end{theorem}
\begin{proof}
The one-variable case is Proposition~\ref{prop:one-dimensional} below. A quadratic polynomial with a local minimum has a positive semidefinite quadratic part, so its minimum is global; see Proposition~\ref{prop:quadratic}. By Theorem~\ref{thm:cubic}, a polynomial of exact degree three in any dimension cannot have a local minimum as its only critical point. For $n=2$ and $d=4$, the affirmative assertion is precisely the quartic result of Calvert and Vamanamurthy \cite[Theorem~5]{CalvertVamanamurthy1980}.

For $n\ge3$ and $d\ge4$, Corollary~\ref{cor:high-dim} gives a counterexample of degree four, hence also of degree at most $d$. Finally, for $n=2$ and $d\ge5$, the counterexamples of \cite[Theorem~1]{CalvertVamanamurthy1980} apply. This proves the classification.
\end{proof}

\begin{remark}\label{rem:theorem5}
The entry $n=2$, $d=4$ is the only one that is not proved in the present paper; it rests on \cite[Theorem~5]{CalvertVamanamurthy1980}. The proof given there proceeds by a case analysis of topological degrees, several of which are left to the reader.
\end{remark}

The threshold can be summarized as follows. The entry ``Yes'' means that the assertion holds for every polynomial of degree at most $d$ satisfying \eqref{eq:question}; ``No'' means that a counterexample exists.
\begin{table}[ht]
\centering
\caption{Validity of the local-to-global implication by dimension and degree.}\label{tab:classification}
\begin{tabular}{lcccc}
\toprule
& $d\le2$ & $d=3$ & $d=4$ & $d\ge5$\\
\midrule
$n=1$ & Yes & Yes & Yes & Yes\\
$n=2$ & Yes & Yes & Yes & No\\
$n\ge3$ & Yes & Yes & No & No\\
\bottomrule
\end{tabular}
\end{table}

\begin{remark}\label{rem:cubic-vacuity}
For a polynomial of \emph{exact} degree three, the premise of $\mathcal P(n,3)$ cannot hold: a local minimum cannot be its sole critical point (Theorem~\ref{thm:cubic}). Thus the affirmative cubic entries in Table~\ref{tab:classification} should not be interpreted as assertions about the existence of cubics with a unique critical point and a local minimum.
\end{remark}

\section{Elementary proofs for degrees one and two}\label{sec:elementary}
For completeness, we record short proofs that do not require topological degree or an index calculation.

\begin{proposition}\label{prop:one-dimensional}
Let $f:\R\to\R$ be of class $C^1$. If $0$ is its only critical point and a local minimum, then $0$ is a global minimum.
\end{proposition}
\begin{proof}
Since $f'$ is continuous and never vanishes on either interval $(-\infty,0)$ or $(0,\infty)$, it has a constant sign on each. The local minimum forces $f'<0$ on $(-\infty,0)$ and $f'>0$ on $(0,\infty)$; either opposite sign would contradict local minimality. Thus $f$ decreases up to the origin and increases afterwards, giving $f(t)\ge f(0)$ for every $t\in\R$.
\end{proof}

\begin{proposition}\label{prop:quadratic}
If a quadratic polynomial $p:\R^n\to\R$ has a local minimum at a critical point, that point is a global minimum.
\end{proposition}
\begin{proof}
After translating the critical point to the origin and subtracting $p(0)$, the linear term vanishes, and $p(x)=Q(x)$ for a quadratic form $Q$. Local minimality implies $Q(x)\ge0$ for every $x$: otherwise $Q(tx)=t^2Q(x)<0$ for arbitrarily small $t\ne0$. Hence the minimum is global.
\end{proof}

\section{Relation to the original degree-five construction}\label{sec:construction}
The simplest two-variable counterexample from \cite[Theorem~1]{CalvertVamanamurthy1980} is
\begin{equation}\label{eq:quintic}
f(x,y)=y^2+x^2(1+y)^3.
\end{equation}
It has degree five, a positive definite quadratic part at the origin, and no critical point other than the origin. Nevertheless, $f(t,-2)=4-t^2\to-\infty$.

The polynomial \eqref{eq:main-example} is connected to \eqref{eq:quintic} by the identity
\begin{equation}\label{eq:restriction}
p\bigl(x,y,x(1+y)\bigr)=f(x,y).
\end{equation}
Thus the new variable $z$ replaces the quadratic expression $x(1+y)$ in the term $(1+y)z^2$, while the additional nonnegative term $(z-x(1+y))^2$ controls the critical-point equations. Introducing one variable lowers the degree from five to four without losing the strict local minimum or uniqueness of the critical point.

This construction shows that the failure of \eqref{eq:question} in dimension three is not caused by a degenerate minimum: the quadratic part at the unique critical point is positive definite. Combined with the two-variable theorem in \cite{CalvertVamanamurthy1980}, it also establishes that three is the smallest dimension in which a quartic counterexample is possible.

\section*{Declaration of competing interest}
The author declares that he has no known competing financial interests or personal relationships that could have appeared to influence the work reported in this paper.

\section*{Funding}
This research did not receive any specific grant from funding agencies in the public, commercial, or not-for-profit sectors.

\section*{Data availability}
No data was used for the research described in the article.

\section*{Declaration of generative AI and AI-assisted technologies in the manuscript preparation process}
During the preparation of this work the author used Claude Opus 5.5 (Anthropic) to assist with the \LaTeX{} typesetting, to improve the English, and to check the proofs and suggest clarifications, additional remarks and references. After using this tool, the author reviewed and edited the content as needed and takes full responsibility for the content of the published article.


\begin{thebibliography}{7}
\bibitem{RosenholtzSmylie1985}
I. Rosenholtz and L. Smylie,
The ``only critical point in town'' test,
\emph{Math. Mag.} \textbf{58} (1985), no.~3, 149--150.

\bibitem{AshSexton1985}
J. M. Ash and H. Sexton,
A surface with one local minimum,
\emph{Math. Mag.} \textbf{58} (1985), no.~3, 147--149.

\bibitem{Durfee1998}
A. H. Durfee,
The index of $\operatorname{grad} f(x,y)$,
\emph{Topology} \textbf{37} (1998), no.~6, 1339--1361.
\newblock \href{https://doi.org/10.1016/S0040-9383(97)00089-X}{doi:10.1016/S0040-9383(97)00089-X}.

\bibitem{Bisgard2015}
J. Bisgard,
Mountain passes and saddle points,
\emph{SIAM Rev.} \textbf{57} (2015), no.~2, 275--292.
\newblock \href{https://doi.org/10.1137/140963510}{doi:10.1137/140963510}.

\bibitem{Fernando2020}
J. F. Fernando,
On the set of local extrema of a subanalytic function,
\emph{Collect. Math.} \textbf{71} (2020), no.~1, 1--24.
\newblock \href{https://doi.org/10.1007/s13348-019-00245-6}{doi:10.1007/s13348-019-00245-6}.

\bibitem{AhmadiZhang2022}
A. A. Ahmadi and J. Zhang,
Complexity aspects of local minima and related notions,
\emph{Adv. Math.} \textbf{397} (2022), Paper No.~108119.
\newblock \href{https://doi.org/10.1016/j.aim.2021.108119}{doi:10.1016/j.aim.2021.108119}.

\bibitem{CalvertVamanamurthy1980}
B. Calvert and M. K. Vamanamurthy,
Local and global extrema for functions of several variables,
\emph{J. Austral. Math. Soc. Ser. A} \textbf{29} (1980), no.~3, 362--368.
\newblock \href{https://doi.org/10.1017/S1446788700021364}{doi:10.1017/S1446788700021364}.
\end{thebibliography}
\end{document}